\documentclass{notices}
\usepackage{amsfonts,amssymb,amsmath,amscd}
\usepackage{graphicx} 
\usepackage{xcolor}
\usepackage{hyperref}
\usepackage[normalem]{ulem} 

\usepackage{titlesec}
\titleformat{\section}{\sffamily\large\bfseries}{\thesection}{.5em}{}
\titleformat{\subsection}{\bfseries}{\thesubsection}{.5em}{}

\AtBeginEnvironment{biblist}{\catcode`\#=12 }

\title{\emph{Princ-wiki-a Mathematica}: Wikipedia editing and mathematics}

\author{David Eppstein\affil{David Eppstein is a distinguished professor of computer science at the University of California, Irvine.  His e-mail address is david.eppstein@gmail.com. Work of David Eppstein is supported in part by NSF grant CCF-2212129.}
\and
Joel Brewster Lewis\affil{Joel Lewis is an associate professor of mathematics at the George Washington University.  His e-mail address is jblewis@gwu.edu. Work of Joel Lewis is supported in part by Simons Foundation grant 634530.}
\and
Russ Woodroofe\affil{Russ Woodroofe is a professor of mathematics at the University of Primorska.  His e-mail address is russ.woodroofe@famnit.upr.si. Work of Russ Woodroofe is supported in part by the Slovenian Research Agency (research program P1-0285 and research projects N1-0160, J1-2451, J1-3003 and J1-50000).} \\
\and
XOR'easter\affil{XOR'easter is a pseudonymous Wikipedia editor who is active in the Physics and Mathematics WikiProjects.}
}
\date{}

\begin{document}

\maketitle

\section{Introduction}
Over the past 20 years, Wikipedia has gone from a rather outlandish idea to a major reference work, with more than 60 million articles across all languages, including nearly 7 million in English~\cite{size}.  Around 27,000 of these articles concern mathematics~\cite{WPM-banner}, and Wikipedia is the first place that many of us go to learn about a new mathematical idea.  In this overview, we will discuss how to go about creating or editing an article on a mathematical subject.  (Most of this applies equally to topics from other technical fields.)  We will also discuss biographies of mathematicians, articles on mathematical books, and the social dynamics of the Wikipedia editor community.

\section{Who edits mathematics on Wikipedia?}

Wikipedia is owned by a non-profit organization (the Wikimedia Foundation), but the encyclopedia is written entirely by its community of volunteers.  The community is governed by a consensus model that encourages boldly making edits and additions, followed by discussion in the case of disagreements.  Anyone with an internet connection can join this community of editors by making a free account.  Some contributors participate using their real names, others under mild pseudonyms, and others strictly anonymously (including the fourth author of this article).  By most measures, activity on Wikipedia peaked around 2007; for example, more than 60,000 editors made five or more edits on the encyclopedia in April 2007, compared with 35--40,000 each month in 2023 \cite{wikimedia-stats}.

Some editors participate in informal groupings, such as WikiProject Mathematics.  This project includes a discussion board\footnote{\url{https://en.wikipedia.org/wiki/Wikipedia_talk:WikiProject_Mathematics}} where editors can solicit input and discuss individual articles or issues of broader relevance to the mathematical editing community.  This forum saw about 1400 edits in 2023, from about 150 individuals (including unregistered editors)~\cite{WPM-edits}.  Contributors include professional mathematicians and students of mathematics, as well as many people with no special credentials, who nonetheless take an interest in sharing mathematical knowledge with the world at large. It is difficult to be more precise about the editor demographics, due to the pseudonymity mentioned earlier.

\section{A mini-glossary}

One of the authors has written a pamphlet titled ``So, you've decided to write about physics and/or mathematics on Wikipedia''~\cite{XOR-essay}, which includes a brief glossary of Wikipedia jargon. We reproduce that list here with a few modifications. Wikipedians have given some ordinary words new, specialized definitions that are related to but distinct from their everyday meanings. Mathematicians should feel right at home with this! Most of Wikipedia's house rules do make sense, after one thinks about the goals and resources of the project. Even so, the terminology can be intimidating.
\begin{itemize}
    \item \textsc{Reliable sources}: published sources of information that are trustworthy enough for our purposes. Textbooks, journal articles, and monographs from major publishers are generally reliable, but beyond that, deeming a source ``reliable'' can be a context-dependent judgment. For example, a self-published source can be reliable for a non-controversial, non-self-aggrandizing claim about the writer (e.g., ``I joined the Miskatonic University faculty in 2018''), but it is not a source we can use for a claim about someone else.
\item \textsc{No original research}: Wikipedia isn't set up to evaluate claims about new ideas or discoveries. So, we have a rule that we don't say anything unless someone else has published it first. This also means that we don't put together ideas to draw new conclusions: no \textit{synthesis}, even if the component ideas are already out there on the record.  In the context of mathematics, this means that adding a proof of a statement is not a substitute for providing a citation to a reliable source.
\item \textsc{Neutral Point Of View}: taking the viewpoints present in reliable sources and representing them ``fairly, proportionately, and, as far as possible, without editorial bias''. This is \textit{not} ``neutrality'' in the sense of saying one bad thing for every good thing, which would lead to absurdities.  This principle has come into play in the Wikipedia article on Inter-universal Teichm\"uller theory,\footnote{\url{https://en.wikipedia.org/wiki/Inter-universal_Teichm\%C3\%BCller_theory}} for example.
\item \textsc{Notability}: similar to, but more specific than, the everyday idea of ``noteworthiness'', a topic is ``notable'' if it is a suitable subject for a Wikipedia article. This depends upon reliable sources existing on the topic, and on the topic being sufficiently distinct from related ones that it makes sense to have a separate page.  For a mathematical concept, this might include appearing in publications by multiple disjoint groups of authors (at a minimum).
\item \textsc{Due weight}: related to, but distinct from, ``notability''. This refers to the question of how much emphasis to give a point \textit{within} an article.
\item \textsc{Good articles}: articles that have gone through an assessment process and been designated well-written, well-cited, and reasonably comprehensive. \textsc{Featured articles} are a level beyond that, being evaluated by more reviewers at once, and meant to cover their subjects comprehensively and represent the pinnacle of what the encyclopedia can provide. 
\item \textsc{Conflict Of Interest} (COI): an editor having a relationship, financial or otherwise, with the subject matter of an article.
\item \textsc{Talk page}: an auxiliary page associated with an article where modifications to that article can be discussed. There is no hard-and-fast rule for when to jump in and edit an article versus asking about it on the talk page first; generally, the more complicated or potentially controversial a planned revision seems, the better it is to raise a question on the talk page. These pages are for discussion of the \emph{article text}, not the subject the article is about. For example, the talk page for the perfect number article\footnote{\url{https://en.wikipedia.org/wiki/Talk:Perfect_number}} is not a suitable place to speculate about whether any odd perfect numbers exist.  
\end{itemize}
Wikipedia has a \textsc{Manual of style} that covers many issues. Part of this resource is a mathematical style guide\footnote{\url{https://en.wikipedia.org/wiki/MOS:MATH}}. But in many cases, consistency across multiple articles is too much to ask for, and the best one can achieve in practical terms is to have each article be coherent within itself.

\section{The ideal Wikipedia article on mathematics}

In order to discuss the challenges of writing a Wikipedia article about a mathematical topic and making that article worth reading, it will be convenient to have an example that illustrates the most critical problems that an aspiring Wikipedia editor will face. We will accordingly discuss what would be involved in crafting a Wikipedia page about \textit{non-Riemannian hypersquares}, a topic that has all the attributes of a modern mathematical subject except for the minor detail that it is completely fictitious~\cite{ideal-mathematician}.

So: Alice, a postdoc or a junior faculty member, is familiar with non-Riemannian hypersquares and discovers one day that Wikipedia does not have an article about them. She naturally wants to fill the gap. We will assume that Alice has followed the standard advice for Wikipedia editors to start with smaller edits, and through doing so has become familiar with the mechanics of Wikipedia editing.\footnote{For one introduction to these mechanics, see \url{https://blog.mollywhite.net/become-a-wikipedian-transcript/}. A convenient list of existing mathematics articles in need of small or larger edits is automatically updated at \url{https://bambots.brucemyers.com/cwb/bycat/Mathematics.html}.} The first question she must face is \textit{whether the topic is suitable for Wikipedia}. If non-Riemannian hypersquares were just something that her research group cooked up during their weekly meetings, then Wikipedia would not be the place to promote them. Fortunately, they are niche but not that niche, having been an active research area for several decades. In other words, they are a topic Alice knows about, not one she invented. Alice has on her shelf a yellow Springer book about them, from the days when Springer invested in typesetting. There are published articles to draw upon, in journals with solid peer review; those from the past several years are available on the arXiv as preprints. As a nice bonus, there is even a section in a \textit{Reviews of Modern Physics} paper that explains a conjectured application to black-hole thermodynamics. Alice has the materials with which to build a Wikipedian treatment. But how to go about it?

Wikipedia is an encyclopedia, not a textbook. This is a distinction that manifests in big ways and small. For one thing, Wikipedia pages do not come with homework problems at the end. At the level of word choice, the academic \textit{we} is downplayed, as are any flourishes of opinion. It is not encyclopedic to say, ``The most elegant construction was provided by Ayanami''; at most, one could write, ``The construction by Ayanami was praised for clarifying the subject'' with citations to passages that give such praise. Even simple observations should be justified by citations, not by proofs. ``We can easily see'' is out for several reasons. Alice must break habits that years of reading and writing mathematical prose have instilled.

Beyond issues of prose style, the organization of a Wikipedia article is not like that of an academic paper, a textbook chapter, or any other kind of scholarly writing. The primary function of a journal article is to report results, most likely divorced from the original way in which they were found, and taking much of the motivation for granted. (Consider Elric's proof that all non-Riemannian hypersquares of even order are symplectic. It is well known that when Elric first presented this result at an AMS sectional meeting, the proof relied upon advanced techniques in mimetic topology, and the version eventually published in the \textit{Inventiones} was much simplified.) History and motivation, although omitted from research papers, are suitable material for a Wikipedia page, provided that they can be documented.

Wikipedia articles on technical topics are generally meant to proceed from easier, more widely accessible material to more difficult and demanding. This is not necessarily the same as covering a subject in historical order. The early work on non-Riemannian hypersquares grew out of questions raised by Noether and Coxeter, but Alice was only able to appreciate that after she had spent quite some time trying to master the theory herself. The first forays into the regularity conjecture, by Steinberg and Bergstein, seemed to contradict each other until they were synthesized by Sha'arawi. Alice has read the Dover reprint of \textit{Sources for Non-Riemannian Hypersquares}, and she can tell you it wasn't easy going. 

Above all, the introductory part of a Wikipedia page before the first section break, variously called the \textit{lead} or \textit{lede} or just the \textit{intro}, is meant to be as widely comprehensible as possible. It is expected to be shorter than the typical introduction to a mathematics paper, while serving the purposes of both the abstract and the introduction. The goal of providing for both general and specialized audiences can be hard to achieve. For example, there is no popularization of non-Riemannian hypersquares: nothing in the books of Gardner or Cheng, no videos by Vi Hart, and \textit{Quanta} hasn't even tried. There is no published source for a simplifying analogy. An article on topology can rely on the joke about how a topologist can't tell the difference between a coffee cup and a doughnut, but no reliable source has yet made a joke about non-Riemannian hypersquare specialists. (It turns out that they can't tell a bike lock from a Linzer torte.)

Facing all these challenges, it is understandable that Alice's first attempt is a rather short page. If her colleagues see it, they find clues enough to let them situate the topic with respect to what they already know. The article is brief, reliant upon jargon, and lacking a gentle introduction, but it probably does provide what they are looking for. The students and hobbyists, in contrast, are out of luck. They see ``In mathematics,'' the topic name in bold, and then a string of incantations. Understandably but unfortunately, experts tend to have a more benign view regarding the state of Wikipedia than the students do.

What if Alice wishes to be recognized for her Wikipedian achievements? First, she can try earning markers of quality for her article within the Wikipedia community itself. She can push her article to Good or even Featured status.\footnote{Some examples of Wikipedia articles that have reached one of these statuses include 
\href{https://en.wikipedia.org/wiki/Affine_symmetric_group}{Affine symmetric group},
\href{https://en.wikipedia.org/wiki/Random_binary_tree}{Random binary tree}, and
\href{https://en.wikipedia.org/wiki/Speed_of_light}{Speed of light}.} Wikipedia's behind-the-scenes pages explain the criteria in detail, but the greatest challenge that Alice faces in going for either status is not spelled out: \textit{there just aren't many Wikipedia editors with the background knowledge} necessary to referee an article about non-Riemannian hypersquares. Feedback is slow to arrive, and if Alice is unlucky, it is frustratingly shallow when it does. (``Provide a picture? Just how many dimensions do they think this hypersquare lives in?!'') Alice soon finds that writing for Wikipedia requires juggling multiple audiences in a way that most mathematical prose does not.

\begin{figure*}
    \centering
    \includegraphics[width=4.75in]{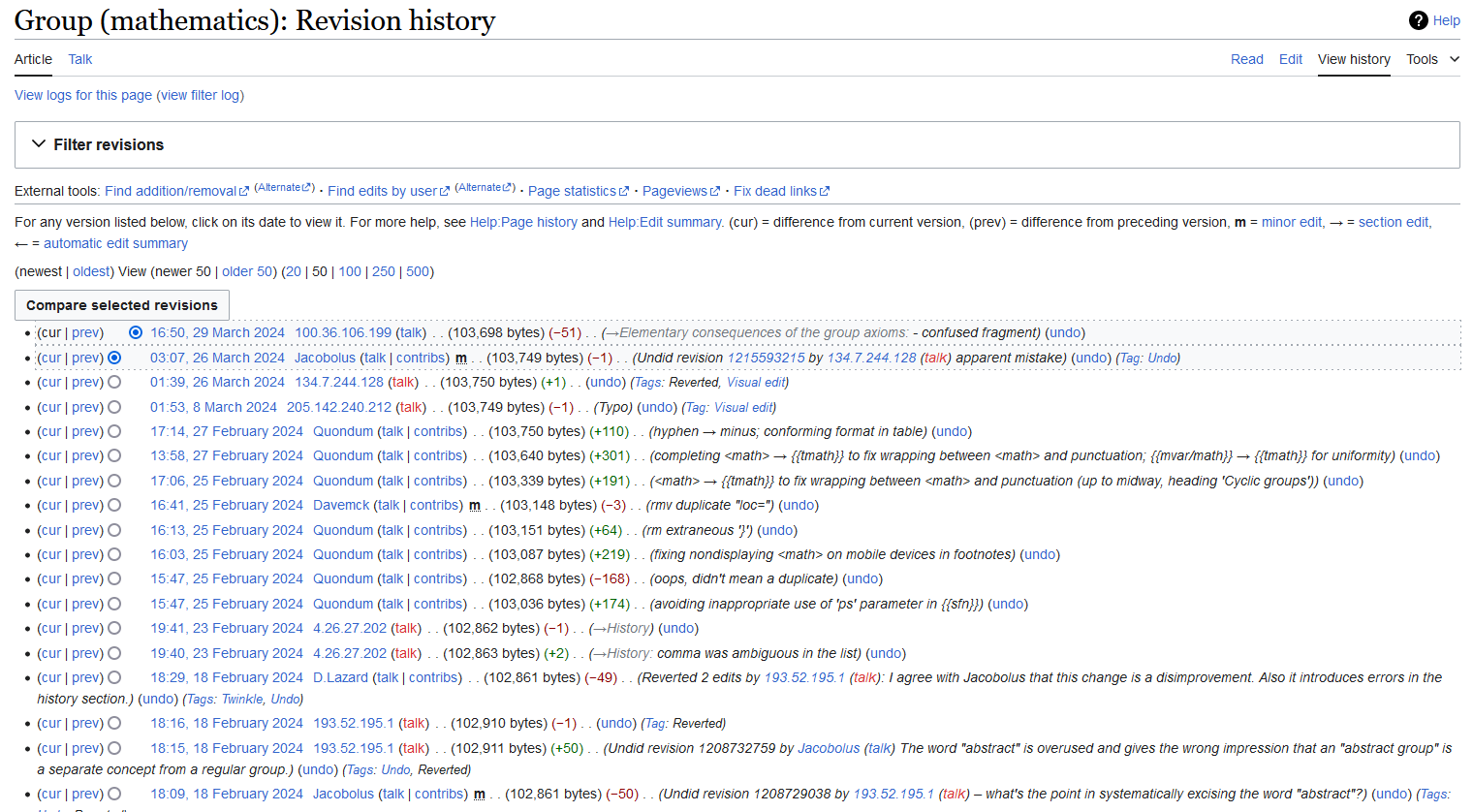}
    \caption{Every article on Wikipedia has all its past versions stored in a viewable record, its revision history}
\end{figure*}

What about getting recognition elsewhere for her contributions to Wikipedia? After all, the tenure committee is waiting, and they've made noises about wanting to see ``outreach'' activities. Wikipedia's design makes this difficult, since the records of who wrote what are buried in the page history. It is also possible, indeed advisable, that Alice edits under a pseudonym. She is a woman on the Internet, after all, with everything that implies. Moreover, she knows that she might get involved in deleting someone's claim to have proved the Collatz conjecture or to have invented a quantum negentropic space drive, or in trimming puffery from an overly promotional autobiography, and the last thing her department chair needs is angry letters from disgruntled parties (a thing that has happened to more than one of the present authors).

\section{Non-research mathematics}

\begin{figure*}
    \centering
    \includegraphics[width=4.75in]{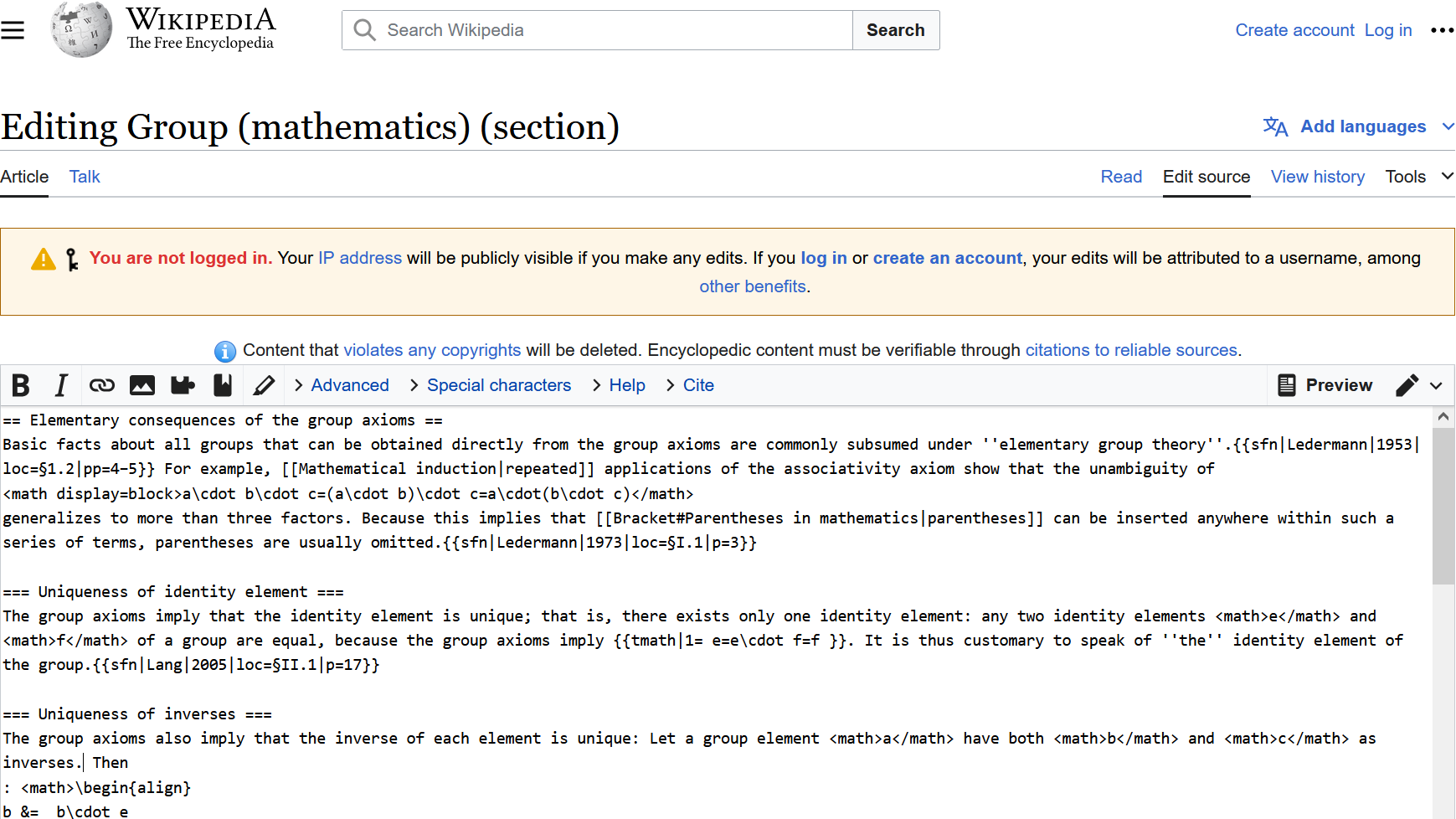}
    \caption{One version of the Wikipedia editing interface, including a bespoke LaTeX implementation}
\end{figure*}

In the previous section, we discussed writing Wikipedia articles about esoteric subjects.  Not all mathematics topics are like that. The Wikipedia article on the number 0,\footnote{\url{https://en.wikipedia.org/wiki/0}} for example, routinely exceeds 200,000 views per month. The article on $\pi$ does pretty well too,\footnote{\url{https://en.wikipedia.org/wiki/Pi}} pulling in around 150,000 views per month, with spikes to more than twice that every March. Wikipedia articles on mathematical topics of more widespread interest such as these present their own challenges. There is probably no topic that is \textit{exclusively} studied at the primary-school level. One can explain prime numbers to a schoolchild---the sieve of Eratosthenes is, in principle, within the grasp of anyone who has mastered their times table---but an encyclopedia article about prime numbers is going to have to cover much more than the definition and the easiest examples. The concept of ``circle'' is taught at such a young age that the reader may not even remember a time before they had learned of it, but an encyclopedic treatment must span Book III of Euclid's \textit{Elements}, coordinate geometry, trigonometry, and even what $S^1$ means to a topologist.

How, then, can an article be all things to all people? What do we do when the same page will be read by PhD students trying to remember a detail and by pre-teens desperate for homework help? This is where the general advice to put the least technical portions up front comes into play. That can, however, collide with the Wikipedia house style of having the lead paragraphs provide a capsule summary of the whole article that follows. The difficulty of the resulting text can escalate towards the end of the lead, then drop back in the first section, and then escalate again as the article proceeds.

The density, placement, and to some extent the purpose of citations differs between Wikipedia and other kinds of mathematical prose. In a textbook, for example, long stretches can go without references because the argument holds together by its own logic. Citations are employed to indicate historical priority and to offload details that the author does not wish to exposit themselves. The question of what topics to include and how long to spend on each one is ultimately up to the writer's discretion. In the Wikipedia realm, however, each claim must be cited to sources that are trustworthy for the purpose for which they are being invoked.  (See ``Reliable sources" in the mini-glossary.)  Moreover, deciding what to include in the first place and how much to write about it also depends upon the literature. Many things can be said about circles; few of them truly need to be inserted into an encyclopedia article called \textsc{Circle}. Accordingly, an observation that might pass un-footnoted in a textbook, or be given with a short proof, might well require one or more citations on Wikipedia, not just to certify its correctness, but to indicate why it was included. For these reasons, well-written Wikipedia articles typically have at least one footnote per paragraph, and often at least one per sentence. 

\section{Biographies of mathematicians}

It's easy to write about famous, dead mathematicians.  There are published obituaries or biographies to use as sources.  The biggest difficulty is cutting down the material to an encyclopedic length. But in these cases, someone has probably already written a Wikipedia article.\footnote{Some examples of good biographies of 20th-century mathematicians include \url{https://en.wikipedia.org/wiki/Ronald_Graham}, \url{https://en.wikipedia.org/wiki/Andrew_M._Gleason}, and \url{https://en.wikipedia.org/wiki/Katherine_Johnson}.}

It's harder to write about working mathematicians, less famous and often still living.  The Wikipedia rules requiring published sources for all content are taken especially seriously for biographies of living people (in Wikipedia jargon, BLPs). The bare facts of a mathematician's degrees and career milestones can be sourced to a publicly available curriculum vitae or to the Mathematics Genealogy Project, but any evaluation of the significance of their accomplishments requires published sources that are independent of the subject and their employers.  Remember that Wikipedia has a strict No Original Research rule!  In many cases, independent coverage of the work of a living mathematicians simply does not yet exist.

In order to provide neutral articles, Wikipedia must guard against all-too-frequent self-promotion, and also against promotion by overly enthusiastic students of subjects.  All biographies should be written by a Wikipedia editor who is separated enough from the subject to not have a conflict of interest.  Wikipedia editors must also evaluate notability, sometimes without having the subject-specific expertise that a tenure committee would have.  As a result, Wikipedia has developed a complex set of criteria to determine which academic researchers are considered notable, and who must wait for an article until their accomplishments become more visible~\cite{Eppstein2018, XOR2023}. A new article on a person deemed unready may well be deleted. Worse, it may be marred for years with a banner telling the world that this person might not be notable.  It's best to be sure of notability before writing a biography.

Some of these criteria are easy to evaluate: Wikipedia should (but doesn't) have articles on all fellows of major academic societies such as the AMS and SIAM, the editors-in-chief of top-tier academic journals, and distinguished or named professors at major research universities. It does have articles on all Fields medalists, but some other major international prize-winners are missing. Authors of multiple notable books (see below) may themselves be notable. In mathematics-adjacent fields such as physics or computer science, where citations are plentiful, notability is often judged by the $h$-index of a researcher or by closely related bibliometrics (such as having many publications with triple-digit citation counts).  Metrics like the $h$-index tend to be less helpful in pure mathematics, where even very significant publications can have few citations, and where major mathematics organizations have strongly discouraged the use of such evaluations~\cite{cite-stats}. If Alice wishes to jump in and create a biography of Beth, it would be wise to have a clear explanation of how she meets the notability criteria.

Wikipedia is not intended as a directory of all professional mathematicians, and being subject of a Wikipedia article is not a reward for having a successful career as a mathematician. Administrative or committee work counts little toward Wikipedia notability, unless it is at the level of the head of an entire university or an entire academic society.  
Associate professors do not usually meet Wikipedia's notability standards, and postdoctoral researchers or assistant professors almost never do, unless they have been recognized by the sort of international award that makes their rising-star status obvious even to non-mathematicians.  In addition, Wikipedia strongly discourages autobiographies.

Overall, Wikipedia biographies tilt very heavily male, in part due to the heavy representation of professional athletes.  The Women in Red project,\footnote{\url{https://en.wikipedia.org/wiki/Wikipedia:WikiProject_Women_in_Red}} named for the fact that links to nonexistent articles appear by default red instead of blue, aims to address this by identifying non-male persons who meet the notability standards for biographies and then writing articles on them.

\subsection{An ideal Wikipedia biography?}
After finishing her article on non-Riemannian hypersquares, Alice decides to also write a biography of Elric, whose work on non-Riemannian hypersquares she is well familiar with.  Elric is a fellow of the AMS, so he is notable.  She starts a draft page at Draft:Elric\_Kova\v{c}, and works on it over a few days.  Alice has met Elric several times at conferences, so after reading Wikipedia's strict COI guidelines, she declares her conflict and submits her draft to the Articles for Creation process.\footnote{\url{https://en.wikipedia.org/wiki/Wikipedia:Articles_for_creation}}  After a few weeks of waiting, an editor with no connection reviews the page.  After a few tweaks, they move the draft to article space at ``Elric Kova\v{c} (mathematician)", since there is already an article on a YouTuber with the same name. Alice's biography of Elric appears as a Wikipedia article.  While she is waiting, she realizes that her colleague Beth has just become a named chair in their department, in addition to being a highly-cited researcher.  Alice adds Beth to the Women in Red project's sublist of possibly-notable women in mathematics, with a link to her homepage and a comment that she is a named chair.  Shortly after, an editor who watches the sublist (possibly one of the present authors!) creates an article on Beth.

\section{Articles about mathematics books}

A decent heuristic is that if a math book has received two or more substantial reviews in reliable venues---so, more than a blurb in ``books received'', and published somewhere with editorial oversight---then it is probably a suitable topic for an article. These reviews serve a dual purpose: they provide the sources on which any article must be based, and their existence provides evidence that the subject is notable. Because of this second purpose, \emph{MathSciNet} and \emph{zbMATH} tend to count less towards notability than other reviews (like those published by the MAA): they may have the same depth of coverage, but since they are perceived as covering everything, a \emph{MathSciNet} or \emph{zbMATH} review doesn't make a book stand out as significant.

We qualify this rule because there will inevitably be edge cases, like books in a series that are better covered together rather than individually. As Wikipedia is a volunteer project, with many things to write about, many books that do meet this standard do not have articles. Writing about writing overlaps with biography. If a person has written a single article-worthy book and has no other record to speak of, it can make more sense to have an article about the book than about the author. On the other hand,
someone that has written several notable books generally qualifies for a biographical article.

At a minimum, a Wikipedia article about a mathematics volume ought to explain who wrote it, when and how it was published, what topics it covers, what audience it is aimed at, and whether and why its reviewers praised or panned it. Alice might choose to write about a book because she found it passionately inspirational, but encyclopedic fair play means reporting when there is a difference of opinions, and summarizing all sides rather than taking one.


\section{Writing about your own work}

Due to the risk of COI, it is generally a bad idea to write about your own research on Wikipedia. Of course, it is tempting to get ahead in the academic rat race by any means possible, and each of us naturally feels that our work \emph{deserves} to be covered there---if we didn't think it was important, we wouldn't be researching it! When faced with this temptation, consider how you would feel about anyone \emph{else} puffing themselves up in this manner, and recognize that (to paraphrase \emph{Calvin and Hobbes}) we are all someone else to somebody.

Rather than advertising their own super-specialization, experts can make themselves useful by explaining the prerequisites to understanding it. What articles would a student read in order to understand the background and broader context of your research? Wikipedia articles are all works-in-progress, and we can guarantee that some fraction of those articles need attention.  You won't be \emph{automatically} respected as an expert, but part of expertise is knowing the right references to go to and how to summarize them.  This will be greatly appreciated, and will also help protect important content, since unsourced material may be deleted or cut far back, even if correct.

It is also possible to suggest edits for other contributors to make, in venues like the WikiProject Mathematics discussion page: ``I don't want to cite my own work, but I covered the quasithick case of non-Riemannian hypersquares in this paper\ldots.'' This breaks the editor's pseudonymity, of course, which one might understandably be reluctant to do.

\begin{figure*}
    \centering
    \includegraphics[width=4.75in]{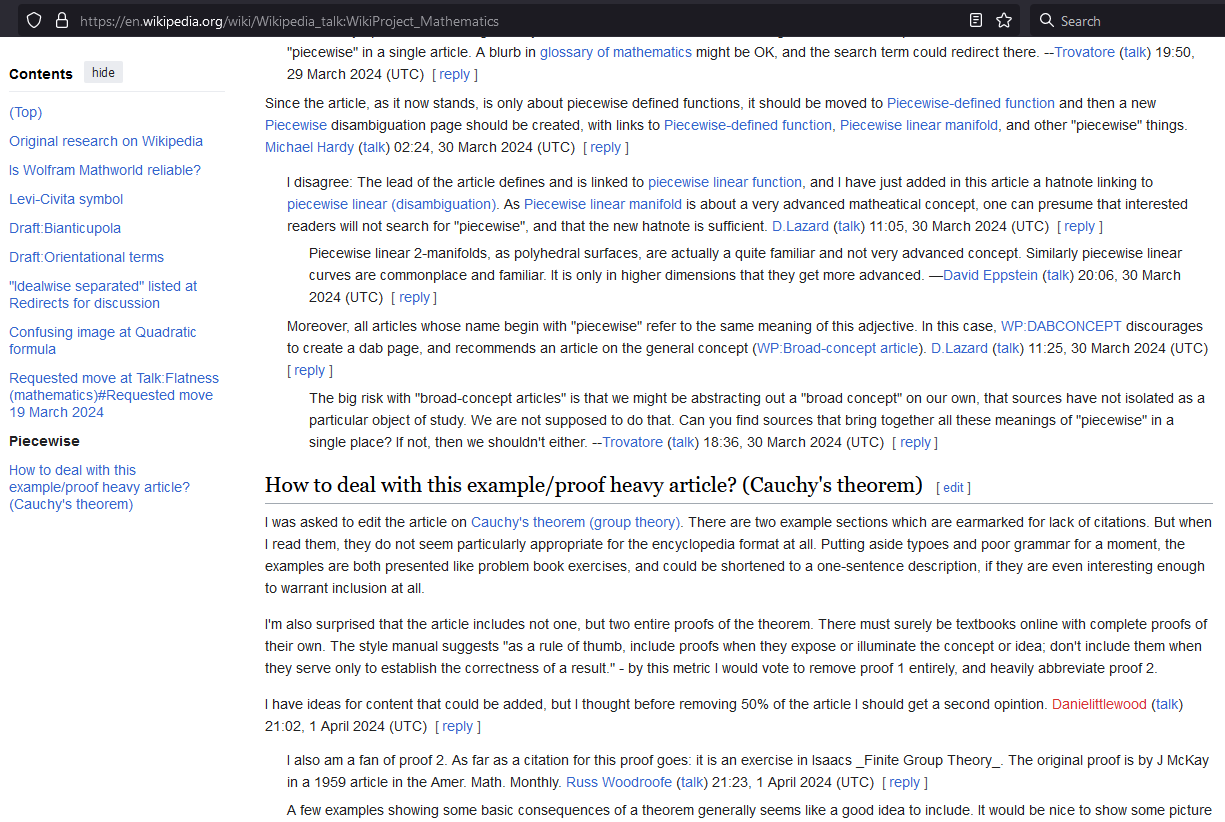}
    \caption{The WikiProject Mathematics page is a centralized discussion point for issues relating to mathematics articles}
\end{figure*}

\section{The community of Wikipedia}

\subsection{Editing fatigue}
The personal satisfaction of editing Wikipedia can be considerable: having one's work appear as the top hit on search engines, read by hundreds of people each day, is exciting in a way that is not easy to duplicate.  Contributing to the most comprehensive reference-work in history is a form of outreach and education on a scale far beyond what we can achieve in the classroom.  It is inspiring to join a community of like-minded people who share in the goal of broadening access to knowledge. In a world drowning in dis- and misinformation, ``I wrote Wikipedia's explanation of calculus'' is what you say just before sliding on your mirrored shades and walking away from an explosion in slow motion.  Indeed, editing can be a little bit addictive, especially to the detail-oriented people that are often drawn to mathematics.

Of course, editing Wikipedia is not always smooth.  In addition to the talk pages associated with articles, there are also some central talk pages.  People who make a serious hobby of editing Wikipedia tend to get into arguments at one or more of these locations. This has all the satisfaction of debating which is the worst \textit{Star Wars} movie or the best comma-placement convention, with the added emotional weight that one is contributing to a reference that people actually use, so it feels like it \textit{matters}. This can lead to ``bikeshedding''~\cite{bikeshedding} over guidelines: lengthy sparring over exactly how to phrase the third bullet point in some documentation page, or exactly how to describe the relationship between one guideline and another. It is easy to tire of this, or, worse, to find oneself enjoying it.  

The need to step away for a time after a contentious talk page argument is common enough that editors talk about the need for a ``wikibreak''. Other strategies include changing area of activity on Wikipedia---there are many articles that need attention. The authors of this article have experienced Wikipedia burnout, occasionally severely so. Managing one's energy and attention is as important on Wikipedia as elsewhere in professional life.

\subsection{Where do editors come from?}
Aside from frustration over on-Wikipedia arguments, another problem is that it can feel like there is too much work to be done in writing an encyclopedia. There is a lot of mathematics, and not enough people to write about it!

Many hands make light work. On the other hand, what ought to be a pleasant hobby can feel like a burden when there are too few people to share it. We have already mentioned the difficulty in obtaining on-wiki feedback about mathematics writing. More generally, it can sometimes feel like a never-ending task: fix one article, and another turns up broken. It is important to remember that you are not responsible for holding up the quality of the project, and easier to remember this if we know that when we cannot be on the case, someone else will be.

The authors do not have a solution. As mentioned previously, early-career mathematicians may have trouble getting credit for their work on Wikipedia. The two of us who edit under our real names came to editing after we had secure jobs. Even then, pseudonymous editing is thankless and editing under your own name may be a mixed blessing. Wikipedia editing might be an excellent retirement project for an emeritus mathematician, and the authors know several emeritus mathematicians who edit somewhat prolifically.

One strategy that we wish to caution against is assigning students the task of editing Wikipedia, at least until you've been a Wikipedia editor yourself for long enough to know what can go wrong. Over the years, many schools have tried involving their students with Wikipedia. This has had mixed results, and many long-time Wikipedia editors have negative reactions to the idea thanks to its being most noticeable when it goes badly. A student with a shaky grasp of the material, writing out of obligation, without the experience to tell a good source from a bad one, and bringing all the habits acquired from essays and term papers \ldots we shudder a little just to think of it. Any instructor contemplating a course that involves editing Wikipedia should, at least, not assign grades based on whether students’ contributions persist: it’s unfair to evaluate students using factors beyond their control. If you are considering involving students in editing, we encourage you to seek the assistance of the WikiEd project.\footnote{\url{https://dashboard.wikiedu.org/}}

\section{Conclusions}

Wikipedia was founded on the outlandish idea that an open, self-selecting community of non-specialists would be able to assemble a meaningful, encyclopedic record of human knowledge.  To an remarkable degree, this idea has turned out to be correct, although not without many complications, difficulties, and flaws.  Wikipedia is now an invaluable and unmatched reference work on mathematics and mathematicians; but one that can still be improved.  The Wikipedia editing community has developed a complicated and always-evolving set of norms over its two decades, and has a complicated relationship to subject matter expertise.  Contributing to Wikipedia presents the opportunity to have your writing reach a much larger audience than is usual for an academic, but under a very different style and model of collaboration.  Wikipedia editing may not be for everyone, but the encyclopedia needs more expert contributors---we hope you will give it a try!



\bibliographystyle{amsplainurl}
\bibliography{refs}

\end{document}